\documentclass[11pt]{article}
\usepackage{amsmath}
\usepackage{amsmath,amssymb,latexsym,color}
\usepackage[mathscr]{eucal}
\newtheorem{thm}{Theorem}[section]

\newtheorem{lemma}[thm]{Lemma}

\newtheorem{example}{Example}[section]
\newtheorem{defin}{Definition}[section]

\newtheorem{remark}{Remark}[section]

\newcommand{\proof}{{\it Proof.\quad}}
\newcommand{\qed}{\hfill\Box\medskip}

\usepackage{CJK}
\begin{document}
\begin{CJK*}{GBK}{song}

\newcommand{\be}{\begin{equation}\label}
\newcommand{\ee}{\end{equation}}
\newcommand{\bea}{\begin{eqnarray}\label}
\newcommand{\eea}{\end{eqnarray}}

\title{\bf The  smallest 3-uniform bi-hypergraphs which are one-realization of a given set }

\author{
Ping Zhao$^{\rm a}$\quad Kefeng Diao$^{\rm a}$\quad Renying
Chang$^{\rm a}$\quad    Kaishun Wang$^{\rm b}$\thanks{Corresponding
author: wangks@bnu.edu.cn}\\
{\footnotesize a. \em School of Science, Linyi  University,
Linyi, Shandong, 276005, China }\\
{\footnotesize b. \em  Sch. Math. Sci. {\rm \&} Lab. Math. Com.
Sys., Beijing Normal University, Beijing 100875,  China} }

\date{}
 \maketitle

\begin{abstract}

For any set $S$ of positive integers, a mixed hypergraph ${\cal H}$
is a one-realization of $S$ if its feasible set is $S$ and each
entry of its chromatic spectrum is either 0 or 1. In this paper, we
determine the minimum size of  3-uniform bi-hypergraphs which are
one-realizations of a given set $S$. As a result, we  partially
solve an open problem proposed by Bujt$\acute{\rm a}$s and Tuza  in
2008.

\end{abstract}

\section{Introduction}

 A {\em mixed hypergraph } on a finite set $X$ is a triple ${\cal
H}=(X, {\cal C}, {\cal D})$, where ${\cal C}$ and ${\cal D}$ are
families of subsets of $X$. The members of ${\cal C}$ and ${\cal D}$
are called {\em ${\cal C}$-edges} and {\em ${\cal D}$-edges},
respectively. A set $B\in {\cal C}\cap {\cal D}$ is called a
\emph{bi-edge}. A \emph{bi-hypergraph} is a mixed hypergraph with
${\cal C}={\cal D}$, denoted by ${\cal H}=(X, {\cal B})$, where
${\cal B}={\cal C}={\cal D}$.
 If ${\cal C}'=\{C\in {\cal C}| C\subseteq X'\}$
 and ${\cal D}'=\{D\in {\cal D}| D\subseteq X'\}$, then the hypergraph ${\cal H}'=(X', {\cal C}', {\cal D}')$  is called a {\em
derived sub-hypergraph}
  of  ${\cal H}$ on $X'$, denoted by ${\cal H}[X']$.

 The distinction between ${\cal C}$-edges and ${\cal D}$-edges becomes substantial when colorings
 are considered. A {\em proper $k$-coloring} of
${\cal H}$ is a partition of $X$ into $k$ \emph{color classes} such
that each ${\cal C}$-edge has two vertices with a {\em Common} color
and each ${\cal D}$-edge has two vertices with {\em Distinct}
colors. A {\em strict $k$-coloring} is a proper $k$-coloring with
$k$ nonempty color classes, and a mixed hypergraph is {\em
$k$-colorable} if it has a strict $k$-coloring. For more
information, see \cite{Tuza, Voloshin1, Voloshin2}.
  The set of all the  values $k$ such that ${\cal H}$ has a strict $k$-coloring is called the {\em feasible set}
  of ${\cal H}$, denoted by $\Phi({\cal H})$. For each $k$, let $r_k$ denote
the number of {\em partitions} of the vertex set. The vector
$R({\cal H})=(r_1,r_2,\ldots,r_{\overline\chi})$ is called the {\em
chromatic spectrum} of ${\cal H}$, where $\overline\chi$ is the
largest possible number of colors in a strict coloring of ${\cal
H}$. If $S$ is a finite set of positive integers, we say that a mixed
hypergraph ${\cal H}$ is a {\it realization} of $S$ if $\Phi({\cal
H})=S$. A mixed hypergraph ${\cal H}$ is a {\it one-realization} of
$S$ if it is a realization of $S$ and all the entries of the
chromatic spectrum of ${\cal H}$ are either 0 or 1. This concept was
firstly introduced by Kr$\acute{\rm a}$l \cite{Kral}.

  It is readily seen that if $1 \in \Phi(\cal{H})$, then $\cal{H}$
cannot have any $\cal{D}$-edges. Let $S$ be a finite set of positive integers with  $\min(S)\geq 2$.
  Jiang et al. \cite{Jiang} proved that the minimum number of
  vertices of realizations  of
$S$  is $2\max(S)-\min(S)$ if $|S|=2$ and $\max(S)-1\notin S$.
Kr$\acute{\rm a}$l \cite{Kral} proved that there exists a 
one-realization of $S$ with at most $|S|+2\max{S}-\min{S}$ vertices.
In \cite{zdw2}, we improved Kr$\acute{\rm a}$l's result and proved
that, for $\min(S)\geq 3$, the smallest size of the one-realizations
of $S$ is $2\max(S)-\min(S)$ if $\max(S)-1\notin
 S$ or $2\max(S)-\min(S)-1$ if $\max(S)-1\in
 S$. Recently, Bujt\'{a}s and Tuza
\cite {Bujtas}
    gave a necessary and sufficient condition for $S$
     being the feasible set of an $r$-uniform mixed
    hypergraph, and they raised the following open problem:

\textbf{Problem.} Determine the minimum number of vertices in
$r$-uniform bi-hypergraphs with a given feasible set.

 In \cite{zdw1}, we constructed a family of 3-uniform bi-hypergraphs with a given feasible set,
 and obtained
an upper bound on the minimum number of vertices of the
one-realizations of a given set.
 In this paper, we focus on this problem and obtain
 the following result:

\begin{thm} For integers $s\geq 2$ and $n_1>n_2>\cdots >n_s\geq 2$,  the smallest size
 of $3$-uniform bi-hypergraphs which are one-realizations of
 $\{n_1,n_2,\ldots,n_s\}$ is
$2n_1-\lfloor \frac{n_2+1}{n_1}\rfloor.$
\end{thm}


\section{Proof of Theorem 1.1}

In this section we always assume that $S=\{n_1,n_2,\ldots,n_s\}$ is
a set of integers with $s\geq 2$ and $n_1>n_2>\cdots >n_s\geq 2$. In  order to prove our main result, we
first give a lower bound on the size of the  3-uniform
bi-hypergraphs which are one-realizations of $S$, then construct two
families of 3-uniform bi-hypergraphs which meet the bound.

\begin{lemma}
 Let $\delta_3(S)$
 denote the minimum size of $3$-uniform bi-hypergraphs ${\cal H}=(X, {\cal B})$  which are one-realizations
 of $S$. Then  $\delta_3(S)\geq 2n_1-\lfloor \frac{n_2+1}{n_1}\rfloor.$
\end{lemma}

\proof We divide our proof into the following two cases.

\textbf{Case 1} $n_1>n_2+1$.

That is to say, $n_1-1\notin S$. Suppose $|X|=2n_1-1$. For any
strict $n_1$-coloring $c=\{C_1,C_2,\ldots,C_{n_1}\}$ of ${\cal H}$,
if there exist two color classes, say $C_1$ and $C_2$, such that
$|C_1|=|C_2|=1$, then $c'=\{C_1\cup C_2,C_3,\ldots,C_{n_1}\}$ is a
strict $(n_1-1)$-coloring of ${\cal H}$, a contradiction. Since
$|X|=2n_1-1$, there exists one color class, say $C_1\in c$, such
that $|C_1|=1$, and $|C_i|=2$ for any $i=2,3,\ldots,n_1$. Suppose
$C_1=\{x\}$ and $C_i=\{x_i,y_i\}, i=2,3,\ldots,n_1$. Then
$c_1=\{\{x,x_2,x_3,\ldots,x_{n_1}\},\{y_2,y_3,\ldots,y_{n_1}\}\}$
and
$c_2=\{\{x_2,x_3,\ldots,x_{n_1}\},\{x,y_2,y_3,\ldots,y_{n_1}\}\}$
are two distinct strict 2-colorings of ${\cal H}$, a contradiction.
If $|X|\leq 2n_1-2$, then we can get a strict $(n_1-1)$-coloring of
${\cal H}$ from a strict $n_1$-coloring of ${\cal H}$, also a
contradiction.

\textbf{Case 2}  $n_1=n_2+1$.

Suppose $|X|=2n_1-2$. For any strict $n_1$-coloring
$c=\{C_1,C_2,\ldots,C_{n_1}\}$ of ${\cal H}$, if there exist three
color classes, say $C_1,C_2$ and $C_3$, such that $|C_1|=|C_2|=|C_3|=1$,
then $c'=\{C_1\cup C_2,C_3,\ldots,C_{n_1}\}$ and $c''=\{C_1, C_2\cup
C_3,C_4,\ldots,C_{n_1}\}$ are two distinct strict $n_2$-colorings of
${\cal H}$, a contradiction. Referring that $|X|=2n_1-2$, there exist two color classes each of which has one
vertex, and each of the other color classes has two vertices. Similar to Case 1,
 ${\cal H}$ has two distinct strict 2-colorings, a contradiction. If $|X|\leq 2n_1-3$, we
 can also have a contradiction.

Hence, the desired result follows. $\qed$

In  the rest,  we shall construct two families of 3-uniform
bi-hypergraphs which meet the bound in Lemma 2.1.

 We first construct the desired 3-uniform
bi-hypergraphs for the case of $n_1>n_2+1$. For any positive integer
$n$, let $[n]$ denote the set $\{1,2,\ldots,n\}$.

\medskip
 {\bf Construction I.}  For $i\in[s]\setminus\{1\},$
write
\begin{eqnarray*}
X_s^s&=&\bigcup\limits_{j=1}^{n_s}
\{(\underbrace{j,j,\ldots,j}_s,0),
(\underbrace{j,j,\ldots,j}_s,1)\},\\
 X_{i-1}^s&=&\bigcup\limits_{k=0}^{n_{i-1}-n_i-1}\{(\underbrace{n_i+k,\ldots,n_i+k}_{i-1},
\underbrace{1,\ldots,1}_{s-i+1},0),
(\underbrace{n_i+k,\ldots,n_i+k}_{i-1},n_i,\ldots,n_s,1)\},\\
 X_{n_1,\ldots,n_s}&=& \{(n_1,n_2,\ldots,n_s,1)\} \cup \bigcup\limits_{t=1}^sX_t^s,\\
{\cal B}_{n_1,\ldots,n_s}&=&\{\{\alpha_1,\alpha_2,\alpha_3\}|
~\alpha_l\in X_{n_1,\ldots,n_s}, l\in [3], |\{\alpha_{1(j)},\alpha_{2(j)},\alpha_{3(j)}\}|=2,j\in [s+1]\}\\
&&\cup \{\{(1,\ldots,1,1,0),(n_s,\ldots,n_s,1,0),
(n_s,\ldots,n_s,n_s,0)\}\},\end{eqnarray*}
 where
$\alpha_{l(j)}$ is the $j$-th entry of the vertex $\alpha_l$.
 Then ${\cal H}_{n_1,\ldots,n_s}=(X_{n_1,\ldots,n_s}, {\cal
B}_{n_1,\ldots,n_s})$ is a 3-uniform bi-hypergraph.

Note that for any $i\in [s]$,
$c_i^s=\{X_{i1}^s,X_{i2}^s,\ldots,X_{in_i}^s\}$ is a strict
$n_i$-coloring of ${\cal H}_{n_1,\ldots,n_s}$, where $X_{ij}^s$
consists of vertices $(x_1,\ldots,x_{i-1},j,x_{i+1},\ldots,x_s,x).$

In the following we shall prove   that $c_1^s,\ldots,c_s^s$ are all
the strict colorings of ${\cal H}_{n_1,\ldots,n_s}$ by induction on
$s$.

\begin{lemma}
 ${\cal H}_{n_1,n_2}$ is a one-realization
of $\{n_1,n_2\}$.
\end{lemma}

\proof Let $c=\{C_1,C_2,\ldots,C_m\}$ be a strict coloring  of
${\cal H}_{n_1,n_2}$. We focuss on the colors of $(1,1,0), (1,1,1)$,
and have the following two possible cases.

\textbf{Case 1} $(1,1,0)$ and $(1,1,1)$ fall into a common color
class.

Suppose $(1,1,0),(1,1,1)\in C_1$. From the bi-edges
$$\{(2,2,0),(1,1,0),(1,1,1)\},
\{(2,2,1),(1,1,1),(1,1,0)\},\{(2,2,0),(2,2,1),(1,1,1)\},$$ we have
$(2,2,0),(2,2,1)\notin C_1$ and $(2,2,0)$ and $(2,2,1)$ fall into a
common color class, say $C_2$. Similarly, we
have $(j,j,0),(j,j,1)\in C_j$ for any $j\in [n_2]$.
The bi-edge $\{(n_2,1,0),(1,1,0),(n_2,n_2,0)\}$ implies that
$(n_2,1,0)\in C_1\cup C_{n_2}$.

\textbf{Case 1.1} $(n_2,1,0)\in C_1$.

The bi-edges $\{(n_2+k,n_2,1),(1,1,1),(1,1,0)\},$
$\{(n_2+k,n_2,1),(n_2,n_2,1),(n_2,1,0)\}$ imply that
$(n_2+k,n_2,1)\in C_{n_2}$ for any $k\in [n_1-n_2]$. Since
$\{(n_2+k,1,0),(n_2+k,n_2,1),(n_2,n_2,1)\},$
$\{(n_2+k,1,0),(n_2+k,n_2,1),(1,1,0)\}$ are bi-edges,
$(n_2+k,1,0)\in C_1$ for any $k\in [n_1-n_2-1]$, and so
$c=c_2^2$.

\textbf{Case 1.2} $(n_2,1,0)\in C_{n_2}$.

For any $j\in [n_2-1]$ and $k\in [n_1-n_2-1]$, the bi-edge
$\{(n_2+k,n_2,1),(j,j,1),(j,j,0)\}$ implies that
$(n_2+k,n_2,1)\notin C_j$; from the bi-edge
$\{(n_2+k,n_2,1),(n_2,1,0),(n_2,n_2,0)\}$, we have
$(n_2+k,n_2,1)\notin C_{n_2}$. Suppose $(n_2+1,n_2,1)\in C_{n_2+1}$.
Since $\{(n_2+1,1,0),(n_2+1,n_2,1),(n_2,1,0)\},
\{(n_2+1,1,0),(n_2,1,0),(n_2,n_2,1)\}$ are bi-edges, $(n_2+1,1,0)\in
C_{n_2+1}$. Similarly, for any $k\in [n_1-n_2-1]$,
$(n_2+k,1,0),(n_2+k,n_2,1)\in C_{n_2+k}$. For any $j\in [n_2-1]$,
the bi-edge $\{(n_1,n_2,1),(j,j,0),(j,j,1)\}$ implies that
$(n_1,n_2,1)\notin C_j$; and for any $k\in [n_1-n_2-1]\cup \{0\}$,
from the bi-edge $\{(n_1,n_2,1),(n_2+k,1,0),(n_2+k,n_2,1)\}$, we
have $(n_1,n_2,1)\notin C_{n_2+k}$. Then $(n_1,n_2,1)\in C_{n_1}$,
and $c=c_1^2$.

\textbf{Case 2} $(1,1,0)$ and $(1,1,1)$ fall into distinct color
classes.

Suppose $(1,1,0)\in C_1,(1,1,1)\in C_2$. From the bi-edge
$\{(n_2,n_2,0),(1,1,0),(1,1,1)\}$, we have $(n_2,n_2,0)\in C_1\cup
C_2$. Suppose $(n_2,n_2,0)\in C_1$. The bi-edges
\begin{eqnarray*}&&\{(n_2,n_2,1),(1,1,1),(1,1,0)\},\{(n_2,n_2,1),(n_2,n_2,0),(1,1,0)\}, \\
 &&\{(n_2,1,0),(n_2,n_2,0),(1,1,1)\}, \{(n_2,1,0),(n_2,n_2,0),(1,1,0)\}
 \end{eqnarray*}imply that $(n_2,n_2,1)\in C_2$ and $(n_2,1,0)\in C_2$. Therefore,
 the three vertices of the bi-edge $\{(n_2,1,0),(n_2,n_2,1),(1,1,1)\}$ fall into a common color class, a contradiction.
 Suppose $(n_2,n_2,0)\in C_2$. Similarly,   we also have a contradiction. It follows that Case 2 does not appear.
$\qed$

\begin{lemma}
 ${\cal H}_{n_1,n_2,n_3}$ is a
one-realization of $\{n_1,n_2,n_3\}$.
\end{lemma}

\proof Let $X_{n_1,n_2,n_3}'=X_3^3\cup X_2^3\cup
\{(n_2,n_2,n_3,1)\}, {\cal H}_{n_1,n_2,n_3}'={\cal
H}_{n_1,n_2,n_3}[X_{n_1,n_2,n_3}']$. Then
$$
\phi:\quad X_{n_1,n_2,n_3}'  \longrightarrow  X_{n_2,n_3},\quad
  (x_2,x_2,x_3,x)  \longmapsto (x_2,x_3,x)
$$
 is an isomorphism from ${\cal H}_{n_1,n_2,n_3}'$ to
${\cal H}_{n_2,n_3}$.

For any strict coloring  $c=\{C_1,C_2,\ldots,C_m\}$  of ${\cal
H}_{n_1,n_2,n_3}$, since the restriction of any strict coloring of
${\cal H}_{n_1,n_2,n_3}$ on $X_{n_1,n_2,n_3}'$ corresponds to a
strict coloring of ${\cal H}_{n_2,n_3}$, by Lemma 2.2 we get the
following two possible cases.

\textbf{Case 1} $(j,j,j,0),(j,j,j,1)\in C_j$ for each $j\in [n_3]$,
 $(n_3+k,n_3+k,1,0)\in C_1$ for any $k\in [n_2-n_3-1]\cup \{0\}$,
and $(n_3+k,n_3+k,n_3,1)\in C_{n_3}$ for any $k\in [n_2-n_3]$.

In this case, we shall prove that $c=c_3^3$. It suffices to discuss
the colors of the vertices in $X_1^3\cup \{(n_1,n_2,n_3,1)\}$.
 The bi-edges $\{(n_2+k,1,1,0),(n_3,n_3,1,0),(n_3,n_3,n_3,1)\},
\{(n_2+k,1,1,0),(n_3,n_3,n_3,0),(n_3,n_3,n_3,1)\}$ imply that
$(n_2+k,1,1,0)\in C_1$ for  any $k\in [n_1-n_2-1]\cup \{0\}$, and from the bi-edges
$\{(n_2+k,n_2,n_3,1),(1,1,1,1),(1,1,1,0)\},
\{(n_2+k,n_2,n_3,1)),(n_3,n_3,n_3,1),(n_3,n_3,1,0)\},$
we have $(n_2+k,n_2,n_3,1)\in C_{n_3}$ for any $k\in [n_1-n_2]$.
Therefore, $c=c_3^3$.

\textbf{Case 2}  $(j,j,j,0),(j,j,j,1)\in C_j$ for each $j\in [n_3]$,
 $(n_3+k,n_3+k,1,0),
(n_3+k,n_3+k,n_3,1)\in C_{n_3+k}$ for any $k\in [n_2-n_3-1]\cup
\{0\}$ and $(n_2,n_2,n_3,1)\in C_{n_2}$.

Then we shall prove that $c=c_1^3$ or  $c=c_2^3$.
Since $\{(n_2,n_2,n_3,1), (n_2,1,1,0), (1,1,1,0)\}$ is a bi-edge, we have
$(n_2,1,1,0) \in C_1\cup C_{n_2}$.

\textbf{Case 2.1} $(n_2,1,1,0)\in C_1$.

For any $k\in [n_1-n_2]$, since $\{(n_2+k,n_2,n_3,1),(n_2,n_2,n_3,1),(n_2,1,1,0)\},
\{(n_2+k,n_2,n_3,1),(1,1,1,0),(1,1,1,1)\}$  are bi-edges, $(n_2+k,n_2,n_3,1)\in C_{n_2}$.
For any $k\in [n_1-n_2-1]$, from the bi-edges $\{(n_2+k,1,1,0),(n_2,1,1,0),(n_2,n_2,n_3,1)\},
\{(n_2+k,1,1,0),
(n_2+k,n_2,n_3,1),(n_2,n_2,n_3,1)\},$ we have
$(n_2+k,1,1,0)\in C_1$. Therefore, $c=c_2^3$.

\textbf{Case 2.2} $(n_2,1,1,0)\in C_{n_2}$.

The bi-edge $\{(n_2+1,n_2,n_3,1),(j,j,j,0),(j,j,j,1)\}$ implies that
$(n_2+1,n_2,n_3,1)\notin C_j$ for any $j\in [n_3-1]$; from the bi-edge
$\{(n_2+1,n_2,n_3,1),(n_3+k,n_3+k,1,0),(n_3+k,n_3+k,n_3,1)\}$, we
have $(n_2+1,n_2,n_3,1)\notin C_{n_3+k}$ for any $k\in [n_2-n_3-1]\cup \{0\}$.
Since $\{(n_2+1,n_2,n_3,1),(n_2,1,1,0),(n_2,n_2,n_3,1)\}$ is a
bi-edge, $(n_2+1,n_2,n_3,1)\notin C_{n_2}$. Suppose
$(n_2+1,n_2,n_3,1)\in C_{n_2+1}$. Then the bi-edges
$\{(n_2+1,1,1,0),(n_2+1,n_2,n_3,1),(n_2,1,1,0)\},
\{(n_2+1,1,1,0),(n_2,1,1,0),(n_2,n_2,n_3,1)\}$ imply
that $(n_2+1,1,1,0)\in C_{n_2+1}$. Similarly, for any $k\in
[n_1-n_2-1]$, $(n_2+k,1,1,0), (n_2+k,n_2,n_3,1)\in C_{n_2+k}$, and
$(n_1,n_2,n_3,1)\in C_{n_1}$. Therefore, $c=c_1^3$.$\qed$

\begin{thm}
${\cal H}_{n_1,\ldots,n_s}$ is a one-realization of
$\{n_1,n_2,\ldots,n_s\}$.
\end{thm}

\proof  By Lemma 2.2 and Lemma 2.3, the conclusion is true for  $s=2$ and $s=3$.
Assume that it is also true for the case of $s-1$.

Let $X_{n_1,\ldots,n_s}'=\bigcup\limits_{i=2}^sX^s_i\cup
\{(n_2,n_2,n_3,\ldots,n_s,1)\}, {\cal H}_{n_1,\ldots,n_s}'={\cal
H}_{n_1,\ldots,n_s}[X_{n_1,\ldots,n_s}']$. Then
$$\psi:\quad  X_{n_1,\ldots,n_s}' \longrightarrow  X_{n_2,n_3,\ldots,n_s},\quad
  (x_2,x_2,x_3,\ldots,x_s,x) \longmapsto~~(x_2,x_3,\ldots,x_s,x)
$$
 is an isomorphism from ${\cal H}_{n_1,\ldots,n_s}'$ to
${\cal H}_{n_2,n_3,\ldots,n_s}$. By induction, all the strict
colorings of ${\cal H}_{n_1,\ldots,n_s}'$ are as follows:
$$c'_i=\{X'_{i1},X'_{i2},\ldots, X'_{in_i}\}, ~~i\in [s]\setminus
\{1\},$$ where $X'_{ij}=X_{n_1,\ldots,n_s}'\cap X_{ij}^s,~~ j\in [n_i]$.

For any  strict coloring $c=\{C_1,C_2,\ldots,C_m\}$  of ${\cal
H}_{n_1,\ldots,n_s}$, since the restriction
 on $X_{n_1,\ldots,n_s}'$  of any strict coloring of
${\cal H}_{n_1,\ldots,n_s}$ corresponds to a strict coloring of
${\cal H}_{n_2,n_3,\ldots,n_s}$, we focus on the restriction of $c$
on  $X_{n_1,\ldots,n_s}'$ and get the following two possible cases:

\textbf{ Case 1} $c|_{X_{n_1,\ldots,n_s}'}=c_2'$.

That is to say $(j,j,x_3,\ldots,x_s,x)\in C_j$ for any $j\in [n_2]$
and $(j,j,x_3,\ldots,x_s,x)\in X_{n_1,\ldots,n_s}'$. In this case,
we shall prove that $c=c_1^s$ or $c=c_2^s$. It suffices  to discuss
the colors of the vertices in $X_1\cup \{(n_1,n_2,\ldots,n_s,1)\}$.
We obtain $(n_2,1,\ldots,1,0)\in C_1\cup C_{n_2}$ from the bi-edge
$\{(n_2,1,\ldots,1,0),(n_2,n_2,n_3,\ldots,n_s,1),(1,\ldots,1,0)\}.$

\textbf{Case 1.1} $(n_2,1,\ldots,1,0)\in C_1$.

For any $k\in [n_1-n_2-1]$, from the bi-edges
\begin{eqnarray*}&&\{(n_2+k,n_2,\ldots,n_s,1),(n_2,n_2,n_3,\ldots,n_s,1),(n_2,1,\ldots,1,0)\},\\
&&\{(n_2+k,n_2,\ldots,n_s,1),(1,1,\ldots,1,1),(1,1,\ldots,1,0)\},\\
&&\{(n_2+k,1,\ldots,1,0),(n_2+k,n_2,\ldots,n_s,1),(n_2,n_2,n_3,\ldots,n_s,1)\},\\
&&\{(n_2+k,1,\ldots,1,0),(n_2+k,n_2,\ldots,n_s,1),(n_2,1,\ldots,1,0)\},\end{eqnarray*} we
have $(n_2+k,n_2,\ldots,n_s,1)\in C_{n_2}$ and $(n_2+k,1,\ldots,1,0)\in C_1$.
The bi-edges \begin{eqnarray*}&&\{(n_1,n_2,\ldots,n_s,1),(n_2,n_2,n_3,\ldots,n_s,1),(n_2,1,\ldots,1,0)\},\\
&&\{(n_1,n_2,\ldots,n_s,1),(1,1,\ldots,1,1),(1,1,\ldots,1,0)\}\end{eqnarray*} imply that
 $(n_1,n_2,\ldots,n_s,1)\in C_{n_2}$. Therefore, $c=c_2^s$.

 \textbf{Case 1.2} $(n_2,1,\ldots,1,0)\in C_{n_2}$.

 For any $j\in [n_s-1]$, the bi-edge
$\{(n_2+1,n_2,\ldots,n_s,1), (j,\ldots,j,0), (j,\ldots,j,1)\}$ implies that $(n_2+1,n_2,\ldots,n_s,1)\notin C_j$.
For any $p\in [s]\setminus \{1,2\}$ and $k\in [n_{p-1}-n_p-1]\cup
\{0\}$, from the bi-edge
$\{(n_2+1,n_2,\ldots,n_s,1),(n_p+k,\ldots,n_p+k,n_p,\ldots,n_s,1), (n_p+k,\ldots,n_p+k,1,\ldots,1,0)\},$ we have
 $(n_2+1,n_2,\ldots,n_s,1)\notin C_{n_p+k}$; and the bi-edge
$\{(n_2+1,n_2,\ldots,n_s,1),
(n_2,n_2,n_3,\ldots,n_s,1),(n_2,1,\ldots,1,0)\}$ implies
$(n_2+1,n_2,\ldots,n_s,1)\notin C_{n_2}$. Suppose
$(n_2+1,n_2,\ldots,n_s,1)\in C_{n_2+1}$. From the bi-edges
\begin{eqnarray*}&&\{(n_2+1,1,\ldots,1,0),(n_2,1,\ldots,1,0),(n_2,n_2,n_3,\ldots,n_s,1)\},\\
&&\{(n_2+1,1,\ldots,1,0),(n_2+1,n_2,n_3,\ldots,n_s,1),(n_2,n_2,n_3,\ldots,n_s,1)\},\end{eqnarray*}
we have $(n_2+1,1,\ldots,1,0)\in C_{n_2+1}$. Similarly, for any
$k\in [n_1-n_2-1]$, $(n_2+k,1,\ldots,1,0),
(n_2+k,n_2,n_3,\ldots,n_s,1)\in C_{n_2+k}$, furthermore,
$(n_1,n_2,n_3,\ldots,n_s,1)\in C_{n_1}$. Therefore, $c=c_1^s$.

\textbf{ Case 2} $c|_{X_{n_1,\ldots,n_s}'}=c_p'$ for some $p\in
[s]\setminus \{1,2\}$.

That is to say,
$(x_2,x_2,x_3,\ldots,x_{p-1},j,x_{p+1},\ldots,x_s,x)\in C_j$ for any
$j\in [n_p]$ and
$(x_2,x_2,x_3,\ldots,x_{p-1},j,x_{p+1},\ldots,x_s,x)\in
X_{n_1,\ldots,n_s}'$. For any $k\in [n_1-n_2]$, we have
$(n_2+k,n_2,n_3,\ldots,n_s,1)\in C_{n_p}$ from the bi-edges
\begin{eqnarray*}&&\{(n_2+k,n_2,\ldots,n_p,\ldots,n_s,1),(1,1,\ldots,1,1),(1,1,\ldots,1,0)\},\\
&&\{(n_2+k,n_2,\ldots,n_p,\ldots,n_s,1),
(n_p,\ldots,n_p,n_{p+1},\ldots,n_s,1),(\underbrace{n_p,\ldots,n_p}_{p-1},1,\ldots,1,0)\}.\end{eqnarray*}
 Then, for any $k\in [n_1-n_2-1]$, from the bi-edges
\begin{eqnarray*}&&\{(n_2+k,1,\ldots,1,0),(n_2+k,n_2,n_3,\ldots,n_s,1),(1,1,\ldots,1,0)\},\\
&&\{(n_2+k,1,\ldots,1,0),(n_2+k,n_2,n_3,\ldots,n_s,1),(n_2,n_2,n_3,\ldots,n_s,1)\},\end{eqnarray*}
we have $(n_2+k,1,\ldots,1,0)\in C_1$. Since
\begin{eqnarray*}&&\{(n_2,1,\ldots,1,0),(n_2,n_2,n_3,\ldots,n_s,1),(1,1,\ldots,1,0)\},\\
&&\{(n_2,1,\ldots,1,0),(n_2,n_2,n_3,\ldots,n_s,1),(n_2+1,n_2,n_3,\ldots,n_s,1)\}\end{eqnarray*} are bi-edges, $(n_2,1,\ldots,1,0)\in C_1$. Therefore, $c=c_p^s$. $\qed$

For the case of $n_2=n_1-1$, we have the following construction.

\medskip
{\bf Construction II.}  Let
$X_{n_1,\ldots,n_s}''=X_{n_1,\ldots,n_s}\setminus
\{(n_2,1,\ldots,1,0)\}$ and ${\cal H}_{n_1,\ldots,n_s}''={\cal
H}_{n_1,\ldots,n_s}[X'']$. Then, for any $i\in [s]$,
$$c_i''=\{X_{i1}'',X_{i2}'',\ldots,X_{in_i}''\}$$ is a strict
$n_i$-coloring of ${\cal H}_{n_1,\ldots,n_s}''$, where
$X_{ij}''=X_{n_1,\ldots,n_s}''\cap X_{ij}^s, j\in [n_i]$.

\begin{thm} ${\cal H}_{n_1,\ldots,n_s}''$ is a one-realization of
$\{n_1,n_2,\ldots,n_s\}$.
\end{thm}

\proof For any strict coloring $c=\{C_1,C_2,\ldots,C_m\}$ of ${\cal H}''$, referring  to the proof of Theorem 2.4,
there are the following two possible cases:

\textbf{ Case 1} $c|_{X_{n_1,\ldots,n_s}'}=c_2'$.

That is to say, $(j,j,x_3,\ldots,x_s,x)\in C_j$ for any $j\in [n_2]$ and
 $(j,j,x_3,\ldots,x_s,x)\in X'$.  Similar to the Case 1 of Theorem 2.4, we have
$(n_1,n_2,\ldots,n_s,1)\in C_{n_2}\cup C_{n_1}$. Therefore,
$c=c_2^s$ if $(n_1,n_2,\ldots,n_s,1)\in C_{n_2}$ and  $c=c_1^s$ if
$(n_1,n_2,\ldots,n_s,1)\in C_{n_1}$.

\textbf{ Case 2} $c|_{X_{n_1,\ldots,n_s}'}=c_p'$ for some $p\in
[s]\setminus \{1,2\}$.

The bi-edges
\begin{eqnarray*}&&\{(n_1,\ldots,n_p,n_{p+1},\ldots,n_s,1),(1,1,\ldots,1,1),(1,1,\ldots,1,0)\}\\
&&\{(n_1,\ldots,n_p,n_{p+1},\ldots,n_s,1),
(n_p,\ldots,n_p,n_{p+1},\ldots,n_s,1),(\underbrace{n_p,\ldots,n_p}_{p-1},1,\ldots,1,0)\}\end{eqnarray*}
imply that $(n_1,n_2,n_3,\ldots,n_s,1)\in C_{n_p}$. Therefore,
$c=c_p^s$. $\qed$

Observe $|X_{n_1,\ldots,n_s}|=2n_1$ and $|X_{n_1,\ldots,n_s}'|=2n_1-1$.
Combining Lemma 2.1, Lemma 2.2, Theorems 2.4 and Theorem 2.5, the proof of
Theorem 1.1 is completed.

\section*{Acknowledgment}

The research
 is supported by NSF of
Shandong Province (No. ZR2009AM013), NCET-08-0052, NSF of China
(10871027) and the Fundamental Research Funds for the Central
Universities of China.

\end{CJK*}


\begin{thebibliography}{99}


\bibitem{Bujtas}
 C. Bujt$\acute{\rm a}$s, Zs. Tuza, Uniform mixed hypergraphs: the possible numbers of colors,
Graphs   Combin. 24 (2008)  1--12.

\bibitem{Jiang} T. Jiang, D. Mubayi, Zs. Tuza, V. Voloshin and D.
West, The chromatic spectrum of mixed hypergraphs,  Graphs  Combin.
18 (2002)  309--318.

\bibitem{Kral}D. Kr$\acute{\rm a}$l, On feasible sets of mixed hypergraphs,
Electron.
J. Combin. 11 (2004)  $\sharp$R19.

\bibitem{Tuza} Zs. Tuza and V. Voloshin, Problems and results on
colorings of mixed hypergraphs, Horizons of Combinatorics, Bolyai
Society Mathematical Studies 17, Springer-Verlag, 2008, pp. 235--255.

\bibitem{Voloshin1} V. Voloshin, On the upper chromatic number of a
hypergraph,   Australas. J. Combin. 11 (1995)  25--45.

\bibitem{Voloshin2} V. Voloshin,  Coloring Mixed Hypergraphs: Theory, Algorithms
and Applications, AMS, Providence, 2002.

\bibitem{zdw1}
P. Zhao, K. Diao and K. Wang, The chromatic spectrum of 3-uniform
bi-hypergraphs, Discrete Math.  311  (2011) 2650--2656.

\bibitem{zdw2}
P. Zhao, K. Diao and K. Wang, The smallest one-realization of a
given set, arxiv: 1106.6099v1  [math. CO].

\end{thebibliography}
\end{document}